\documentclass{article}

\usepackage{arxiv}

\usepackage[utf8]{inputenc} % allow utf-8 input
\usepackage[T1]{fontenc}    % use 8-bit T1 fonts
\usepackage{hyperref}       % hyperlinks
\usepackage{url}            % simple URL typesetting
\usepackage{booktabs}       % professional-quality tables
\usepackage{amsfonts}       % blackboard math symbols
\usepackage{nicefrac}       % compact symbols for 1/2, etc.
\usepackage{microtype}      % microtypography
\usepackage{lipsum}		% Can be removed after putting your text content
\usepackage{graphicx}

\usepackage{moreverb}

\usepackage{graphicx, subfigure}

\newcommand\BibTeX{{\rmfamily B\kern-.05em \textsc{i\kern-.025em b}\kern-.08em
T\kern-.1667em\lower.7ex\hbox{E}\kern-.125emX}}

\usepackage{graphicx}
\RequirePackage{fix-cm}
\usepackage{hyperref}
 \usepackage{epsfig} % for postscript graphics files
 \usepackage{amssymb}
 \usepackage{amsmath}
\usepackage{multirow}
\usepackage{epsfig} % for postscript graphics files
\input{epsf}
\usepackage{graphicx}
\usepackage{graphicx, subfigure}
\usepackage[utf8]{inputenc}
\usepackage[english]{babel}

\def\beq{\begin{equation}}
\def\eeq{\end{equation}}
\def\baq{\begin{eqnarray}}
\def\eaq{\end{eqnarray}}
\def\bal{\begin{align} }
\def\eal{\end{align} }
\def\bc{\begin{center}}
\def\ec{\end{center}}

\def\fai{\varphi}

\def\fai{\varphi}

\def\gr{\gamma_{r}}
\def\gr1{\gamma_{r1}}

\def\fai{\varphi}

\title{RLS Framework with Segmentation of the Forgetting Profile and Low Rank Updates}
\author{ Alexander Stotsky \\
      Independent Researcher \\
      Gothenburg, Sweden \\
        \texttt{alexander.stotsky2@telia.com}  }
\hypersetup{
pdftitle={A template for the arxiv style},
pdfsubject={q-bio.NC, q-bio.QM},
pdfauthor={David S.~Hippocampus, Elias D.~Striatum},
pdfkeywords={First keyword, Second keyword, More},
}
\date{}
\begin{document}
\maketitle

\begin{abstract}
~This report describes a new regularization approach based on segmentation of the forgetting profile in sliding window least squares estimation. Each segment is designed to enforce specific desirable properties of the estimator such as rapidity, desired condition number of the information matrix, accuracy, numerical stability, etc.
The forgetting profile is divided in three segments, where the speed of estimation is ensured by the first segment, which employs rapid exponential forgetting of recent data. The second segment features a decline in the profile and marks the transition to the third segment, characterized by slow exponential forgetting to reduce the condition number of the information matrix using more distant data. Condition number reduction mitigates error propagation, thereby enhancing accuracy and stability.
This approach facilitates the incorporation of a priori information regarding signal characteristics (i.e., the expected behavior of the signal) into the estimator.
\\
Recursive and computationally efficient algorithm with low rank updates based on new matrix inversion lemma for moving window associated with this regularization approach is developed.
\\
New algorithms significantly improve the approximation accuracy of low resolution daily temperature measurements obtained at the Stockholm Old Astronomical Observatory, thereby enhancing the reliability of temperature predictions.
\end{abstract}


\maketitle

\section{Finite Versus Infinite Window: Opportunities \& Challenges}
\label{win}
\noindent
The Recursive Least Squares (RLS) algorithm with exponential forgetting, \cite{fom},\cite{lju1} is often preferred over sliding window least squares in real-time applications due to its recursive structure and lower per-update computational complexity.
RLS algorithm updates estimates with $\mathcal{O}(n^2)$
operations, while moving window least squares requires $\mathcal{O}(n^3)$
complexity due to repeated matrix inversions, making RLS more efficient for online use.
\\
Nevertheless estimation over the finite window with exponential forgetting offers key advantages over classical RLS with infinite memory, as outlined below:
\begin{enumerate}
\item {\it Faster Adaptation to Non-Stationary Environments.}
In practical DSP (Digital Signal Processing) applications, such as channel equalization, echo cancellation, or biomedical DSP, the signals of interest often exhibit non-stationary behavior due to time-varying system characteristics or non-stationary noise environments. Adaptive filters must therefore be capable of tracking these changes effectively. The classical RLS algorithm with exponential forgetting attempts to address non-stationarity by applying a forgetting factor $\lambda < 1$, which exponentially downweights older data. However, this approach still retains an infinite memory: all past samples, even those from the distant past, continue to influence the filter coefficients, though with diminishing weight. This long memory tail can slow down the algorithm's ability to adapt to sudden signal changes or abrupt parameter shifts.
\\
A sliding window variant of the RLS algorithm with exponential forgetting offers a more responsive solution. By enforcing a finite window over the input data,combined with exponential weighting, the algorithm effectively imposes a hard cutoff on memory, completely discarding older observations beyond the window's length. This leads to faster adaptation since only the most recent data, which is more relevant to the current signal state, contributes to the filter update. The result is quicker convergence of filter coefficients and improved tracking performance in rapidly changing environments. Such characteristics make the sliding window RLS particularly advantageous in DSP applications that demand low-latency and real-time responsiveness under non-stationary conditions.

\item {\it Suppression of Long-Term Influence from Outliers and Past Noise.}
In many DSP applications, such as speech enhancement, radar DSP, and real-time sensor data analysis, the presence of outliers, transient disturbances, or bursts of noise can significantly degrade the performance of adaptive filtering algorithms. The classical RLS algorithm with exponential forgetting partially mitigates this by downweighting older observations via a forgetting factor $\lambda < 1$. However, due to its infinite-memory structure, no data point is ever fully discarded. As a result, outliers or high-noise measurements from the distant past continue to exert residual influence on the parameter estimates, potentially leading to biased or sluggish filter behavior, particularly in environments with frequent or unpredictable disturbances.
\\
By contrast, the sliding window RLS framework introduces a finite memory by retaining only the most recent observations within a fixed-length window. This hard truncation ensures that data, including outliers and noise spikes, are completely eliminated from the estimation process after they fall outside the sliding window. As a result, the algorithm becomes inherently more robust to isolated anomalies and better suited for real-time DSP applications where transient noise events are common. Moreover, this approach improves the filter's ability to maintain accurate and stable parameter estimates over time, as outdated or corrupted information is systematically purged. The suppression of long-term influence from such irregularities enhances both the reliability and responsiveness of adaptive systems operating in noisy or unpredictable signal environments.

\item {\it Improved Numerical Stability.}
One of the key challenges in implementing the classical RLS algorithm in practical DSP systems is maintaining numerical stability over long periods of operation. Due to its infinite memory structure, standard RLS continuously incorporates information from all past observations, albeit with exponentially decaying weights. Over time, this indefinite accumulation of data can lead to numerical issues such as the gradual build-up of rounding errors, loss of precision, and degradation in the quality of matrix operations. In particular, the condition number of the  information  matrix may increase significantly, making the recursive matrix inversion more sensitive to numerical errors and potentially leading to unstable behavior or divergence in finite precision implementations.
\\  For example, when running classical RLS  with infinite memory in a low-battery, on-board system, key risks include numerical instability due to unreliable arithmetic operations, data corruption in internal variables (like recursive inversion of the information matrix), and timing issues that disrupt the recursive update cycle. Inaccurate sensor inputs and memory errors can further degrade estimation quality or cause divergence. Since classical RLS is highly sensitive to precision and timely updates, these issues can lead to incorrect parameter tracking or complete algorithm failure if not properly managed.
\\
Sliding window RLS algorithms address this issue by constraining the memory of the system to a finite set of the most recent data points. By discarding older data that no longer contributes meaningfully to current estimation, the algorithm effectively limits the accumulation of numerical errors. This not only prevents the condition number of the information matrix from growing unbounded but also ensures that the matrix remains well-conditioned and invertible over time. The result is a more stable and reliable estimation process, particularly important in real-time DSP applications where long-term algorithmic stability is critical and computations must often be performed in fixed- or limited-precision environments. Thus, the sliding window approach enhances the robustness of RLS implementations by combining fast adaptation with improved numerical behavior.

\item {\it Greater Flexibility in Memory Control.}
In adaptive filtering applications within DSP, such as real-time system identification, noise cancellation, and time-varying channel estimation, the ability to control how much past information influences current parameter estimates is crucial for achieving a balance between responsiveness to changes and stability against noise. The classical RLS algorithm achieves memory control solely through the exponential forgetting factor $\lambda$, which governs the rate at which past data is downweighted. While effective to some degree, this single-parameter mechanism couples memory depth and weighting, limiting the algorithm's flexibility. Specifically, tuning $\lambda$ to improve tracking performance can inadvertently lead to increased estimation noise, and vice versa.
\\
In contrast, the sliding window RLS framework introduces a second, independent mechanism for memory control: the window length. By allowing separate tuning of the forgetting factor and the size of the data window, the algorithm gains a two-tiered structure for managing memory. The window length imposes a hard cutoff on how far into the past the algorithm considers data, directly controlling the estimator’s effective memory span. Meanwhile, the forgetting factor continues to shape the relative importance of recent observations within that window. This decoupling allows for significantly finer control over the trade-off between estimation smoothness and adaptation speed. For example, a short window with moderate forgetting may yield a highly responsive filter, while a longer window with strong forgetting provides smoother estimates with greater noise resilience.
\\
This added flexibility is particularly valuable in dynamic or unpredictable signal environments, where optimal memory settings may vary over time or across applications. It also facilitates more precise algorithm tuning in resource-constrained DSP systems, enabling designers to tailor performance according to system latency, noise characteristics, or computational capacity.

\item {\it Reduction of Bias from Initial Conditions.}
In practical DSP applications, especially those involving real-time operation or deployment in dynamic environments, the initial conditions of adaptive filters can have a significant and sometimes detrimental effect on performance. In classical RLS algorithms with infinite memory, the initialization of parameters, including the initial estimate vector and the inverse of information matrix plays a critical role. If the initial estimates are inaccurate, noisy, or not representative of the actual signal characteristics, their influence can persist for a long time due to the algorithm’s recursive nature and reliance on all past data with exponentially decaying weights. Although the forgetting factor $\lambda < 1$ gradually reduces the impact of earlier data, it does not eliminate it entirely. As a result, these early conditions can introduce a long-lasting bias in the parameter estimates, especially when the signal statistics evolve or change over time.
\\
Sliding window RLS algorithms mitigate this problem by imposing a strict finite-memory structure. Only a limited number of the most recent data samples are retained within the window, while older data, including any influenced by poor initial conditions,  is fully discarded. This means that after a sufficient number of new observations have entered the window, the initial estimates and early errors no longer affect the current parameter updates. As a result, the estimator becomes more robust to initialization choices, and its output more accurately reflects the present characteristics of the signal. This feature is particularly beneficial in scenarios where the system starts from a cold state, is exposed to transient artifacts during startup, or operates in non-stationary environments where early data quickly becomes obsolete.
\\
By effectively resetting the estimator over time, sliding window methods ensure that the adaptive filter remains focused on the most relevant and recent information. This not only reduces the potential for long-term bias but also contributes to faster convergence and improved estimation accuracy in changing signal conditions.

\item {\it Facilitated Theoretical Analysis in Non-Ideal Conditions.}
Theoretical analysis of adaptive filtering algorithms plays a central role in DSP system design, enabling the prediction of convergence behavior, error performance, and robustness under various operating conditions. However, analyzing the classical RLS algorithm, particularly under non-ideal conditions such as non-stationary inputs, time-varying noise, or model mismatch presents significant challenges due to its infinite-memory structure. In standard RLS with exponential forgetting, every past observation continues to exert some influence on the current parameter estimates, making it difficult to isolate the effects of specific data subsets or to define a bounded region of influence for convergence and error analysis. This introduces mathematical complexity in deriving closed form error bounds, transient behavior, or stability guarantees, especially in the presence of noise bursts or structural changes in the signal.
\\
Sliding window RLS methods offer a key advantage in this regard by operating on a finite, well-defined data horizon. The fixed window length imposes a natural bound on the estimation history, which simplifies the underlying mathematical models. With only a limited number of observations influencing the parameter updates at any given time, analytical derivations such as mean-square error bounds, convergence rates, or sensitivity to noise can be made more tractable. The finite-memory structure also allows for more straightforward application of matrix analysis, time-varying system analysis, and robustness evaluations in practical scenarios where ideal assumptions (e.g., stationarity or Gaussianity) do not hold.
\\
Moreover, in scenarios involving abrupt changes, such as switching systems, time-varying filters, or bursty interference the sliding window framework provides a clear analytical boundary after which the effects of outdated or corrupted data are fully eliminated. This allows for tighter performance guarantees and more modular analysis approaches, where the system’s behavior can be studied independently over successive windows.
\\
As a result, sliding window RLS not only offers practical advantages in adaptive DSP systems but also facilitates the development of rigorous theoretical models that better reflect real-world, non-ideal conditions.
\end{enumerate}
The   sliding window approach   has emerged as a compelling framework in the context of adaptive estimation and real-time DSP. Its key strength lies in its ability to provide  localized, time sensitive estimates  by maintaining and operating on a finite, most recent subset of data samples, thereby ensuring adaptability to time-varying environments. This localized processing enables more responsive estimation compared to full-history methods, particularly in non-stationary scenarios.
\\
Recent advancements have significantly enhanced the appeal of this approach. Notably, the development of   computationally efficient recursive algorithms   that operate with   quadratic complexity   has addressed longstanding concerns about the scalability and feasibility of sliding window methods in real-time applications. These algorithms  \cite{sto2023} - \cite{sto2024} maintain computational tractability while still leveraging the statistical advantages offered by windowed estimation, such as improved tracking capabilities and reduced latency in response to changes in the underlying system dynamics.
\\
As a result, the sliding window framework has increasingly been recognized as a   viable and often superior alternative to classical RLS algorithms, especially in applications where computational resources are limited or where rapid adaptation is required. Unlike classical RLS, which relies on the entire historical dataset, sliding window methods inherently discard outdated information, making them more robust to shifts in system behavior and better suited for online learning or streaming data scenarios.
\\
Nevertheless, despite these advantages, the goal of   further improving estimation performance   remains a central focus in the ongoing development of sliding window techniques. Several challenges and limitations   persist, which currently restrict the full potential of these methods. These include issues related to the optimal selection of forgetting factor and  window size, stability and numerical robustness in recursive implementations, sensitivity to noise and model mismatch, and a number of trade-offs.
\\
The following items outline these challenges in greater detail, serving to highlight critical areas where further research and methodological innovations are necessary to unlock the full capabilities of sliding window-based estimation methods in both theoretical and applied contexts:

\begin{enumerate}
\item {\it Rapidity Costs.}
While sliding window and exponentially forgetting RLS algorithms offer significant advantages in terms of rapid adaptation to time varying signals, this speed comes at a cost, both numerically and theoretically. When estimating rapidly changing parameters, the algorithm must rely on a short memory of the data. This is typically achieved through a small window size and a low forgetting factor, which together place a strong emphasis on recent samples. However, this approach introduces several critical trade-offs that affect estimation accuracy, stability, and spectral resolution.
\\
First, the use of a short sliding window can result in a high condition number of the information matrix, due to reduced data diversity and increased collinearity among the regressor vectors. This ill-conditioning amplifies the variance of the parameter estimates and can lead to significant numerical instability, especially in finite-precision implementations. In extreme cases, the matrix may become nearly singular or even rank-deficient, preventing reliable parameter updates and degrading overall filter performance.
\\
Second, from a signal representation perspective, short windows cannot capture full cycles of   low-frequency components, such as slowly varying harmonic content. Over short intervals, sine and cosine regressors with low frequencies become nearly indistinguishable from lower-order polynomials, sine waves may be approximated as linear terms, while cosine components resemble constants. As a result, the estimator naturally   blends harmonic and trend interpretations, allowing it to track fast transients effectively, but at the expense of   true harmonic resolution  for low frequencies. In practice, this means that while the estimator remains responsive to rapid changes, it inherently   filters out or misrepresents  slowly varying harmonics, effectively trading frequency resolution for temporal responsiveness.
\\
This trade-off also introduces   non-uniform variance across parameters, a phenomenon known as   variance imbalance. Parameters corresponding to lower-frequency regressors, which do not complete even a single period within the sliding window tend to have   much higher estimation variances  than those associated with higher-frequency components, whose periodicity is better captured within the same window. Furthermore, when regressors become nearly collinear, the parameter space becomes poorly defined, leading to further inflation of variance and potential divergence in recursive updates.
\\
Despite these limitations, rapid estimators still provide valuable tracking capabilities during transients, often outperforming slower, more stable estimators in scenarios where responsiveness is paramount. However, the   price for rapid tracking is high variance, lost harmonic resolution at low frequencies, and numerical instability  .
\\
Addressing these challenges requires the development of advanced algorithms that preserve fast adaptation while  reducing the condition number, mitigating variance imbalance, and enhancing numerical robustness, particularly in short-window, high-resolution adaptive DSP environments.

\item {\it Reduction of Computational and Memory Costs.}
Sliding window approaches, as opposed to infinite-memory formulations, offer an opportunity to reduce computational and memory costs   in adaptive filtering, especially when designed with efficient recursive updates. In the classical RLS algorithm, each new data sample contributes to the estimation through an accumulation of all past observations, which requires updating the full information matrix recursively. This growing history increases both  computational complexity  and memory requirements, and over time can pose challenges in real-time or embedded DSP implementations.
\\
In contrast, the sliding window RLS algorithm maintains a  fixed-size buffer of recent data, meaning that the estimator only considers a limited number of observations at any given time. This naturally bounds the memory usage, making the method suitable for resource constrained platforms such as embedded processors, real-time DSP systems, and battery-operated devices. However, this structure also introduces the need for efficient computation as the window slides, when   new samples are added  and old samples are removed  at each time step.
\\
Each of these operations, adding a new data point or removing an old one affects the structure of the information matrix. Specifically, the movement of the window involves  updating the matrix by a number of rank-one matrices: one for each outer product between the regressor vectors and their corresponding contributions to the matrix. Efficient algorithms must therefore be developed to  perform recursive updates  to both the information matrix and its inverse, without recomputing them from scratch at every iteration. Tools such as the  Sherman–Morrison–Woodbury identity  and other matrix inversion lemmas are commonly employed to handle these updates involving a number of rank-one modifications.
\\
When designed properly, these recursive methods significantly  reduce computational overhead, enabling real-time processing even in high-dimensional systems. They also  enhance numerical stability, as they avoid large-scale matrix operations and reduce the risk of rounding errors or overflow, which are especially problematic in long-term, high-rate filtering.
\\
Summarizing the facts presented above it is possible to conclude that the sliding window approach adds algorithmic complexity due to the continuous addition and removal of data, but it also facilitates  bounded memory usage, scalable computation, and   adaptability to modern DSP platforms. The key enabler of these benefits is the development of recursive algorithms capable of efficiently managing updates involving   a number of rank-one matrices, ensuring that the system remains responsive and efficient without sacrificing stability or accuracy.
\end{enumerate}

Ill-conditioning of the information matrix, as discussed above, presents a significant challenge in many estimation problems, particularly when dealing with noisy measurements, closely spaced parameters, or ill-posed inverse problems. This numerical instability can lead to unreliable parameter estimates, large variances, or even the failure of the estimation algorithm to converge. To mitigate these issues, various regularization strategies have been developed and widely adopted.

One common approach involves modifying the information matrix by adding small positive values to its diagonal entries, commonly referred to as diagonal loading or Tikhonov regularization, \cite{ker}. This technique improves the condition number of the matrix, thereby enhancing numerical stability during matrix inversion or optimization procedures. However, the choice of regularization parameter (i.e., the magnitude of the diagonal perturbation) plays a crucial role, as overly aggressive regularization can suppress important signal components.
\\
Another effective strategy is to approximate the low-frequency components of the model using polynomial functions. By capturing the dominant trends with a lower-dimensional polynomial basis, this method can reduce the dimensionality of the estimation problem and improve robustness. Similarly, selecting alternative basis functions, such as orthogonal polynomials, wavelets, or splines can provide better numerical conditioning compared to standard basis choices, especially in cases where the original basis leads to nearly linearly dependent columns in the design matrix. These alternatives often exhibit improved orthogonality properties and better behavior under numerical operations.
\\
Other regularization techniques may involve prior information incorporation, sparse representations, or constrained optimization frameworks. Each of these methods introduces a form of bias or constraint into the estimation process, intended to counteract the ill-conditioning.
\\
Importantly, all these regularization approaches involve inherent trade-offs,\cite{bern1}.   While they can significantly enhance numerical stability and convergence properties, they may also degrade estimation accuracy by introducing bias or limiting the model’s flexibility. The degree of performance degradation depends on the nature of the regularization and the specific characteristics of the estimation problem at hand. Therefore, careful design, tuning, and validation of the chosen regularization strategy are essential to balance stability and estimation fidelity.
\begin{figure*}[!ht]
\begin{center}
%\hspace{-2cm} % or -2cm to move left
  \includegraphics[height=10cm,width=15.5cm]{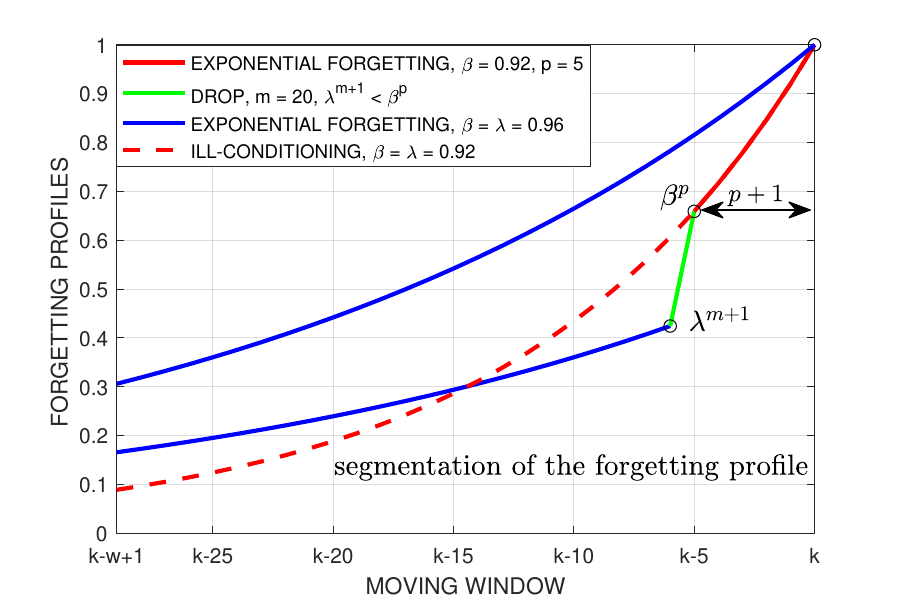}
\end{center}
\caption{\small{The Figure shows the  forgetting profile in the moving window of the size $w$ segmented by red, green and blue lines (RGB lines). The rapidity of estimation is guaranteed by rapid exponential forgetting with small forgetting factor $\beta = 0.92$ associated with the red line (red segment).  Such rapid forgetting over the whole window implies ill-conditioning of the information matrix, see red dashed line.  Reduction of the condition number is associated with the blue line and the forgetting factor $\lambda = 0.96$.  The resulting  profile is segmented by the red segment, drop (with the magnitude determined by a positive   integer $m$) which is plotted with the green line and blue segment with larger forgetting factor. Past data (associated with blue segment) forms the basis for reduction of the condition number of the  information matrix and the variances of parameter estimates, while rapidity is achieved by fast forgetting of recent measurements.}}
\label{figseg1}
\end{figure*}
% everything is in programme ro_6_2_profil2_arxiv1
\section{Segmentation of the Forgetting Profile: A Win-Win Solution for Rapidity and Condition Number}
\vspace{-0.1cm}
\label{fp}
\noindent
The most physically grounded and effective regularization method involves design/segmentation of the forgetting profile in the sliding window, where each segment assigns desired property to the estimator, \cite{slok}. Rapidity, desired condition number of the information matrix, accuracy and numerical stability can be mentioned among desired properties.
The proposed segmentation framework facilitates optimization of the trade-offs associated with the properties of the estimator
and supports the incorporation of prior knowledge about the anticipated properties of the signal into the estimator.
% Rapid forgetting is necessary, for example, in scenarios where abrupt changes in the signal are anticipated.
\\ New method of segmentation of the forgetting profile illustrated in Figure~\ref{figseg1} is described in this report. The rapidity of estimation is guaranteed by fast exponential forgetting, see the red segment (red solid line) in  Figure~\ref{figseg1}. Such rapid forgetting over the whole window (see red dashed line) implies ill-conditioning of the information matrix. Reduction of the condition number is associated with the blue line, which represents slow exponential forgetting. Rapid estimation can not be achieved with blue forgetting profile, but this profile mitigates ill-conditioning and improves robustness and accuracy.
The segmented profile introduces the drop in the red line (which can be adjusted), followed by transition to the blue segment/tail, where the profile diminishes more slowly due to a higher forgetting factor.
The condition number of the information matrix is reduced using more distant data (due to the blue segment) which implies lower variance of parameter estimates and improved numerical stability, while the rapidity is achieved by fast forgetting of recent measurements (due to the red segment).

\section{Matrix Inversion Lemma for Moving Window}
%
%\section{Recursive Estimation Framework for Segmented Forgetting Profile}
%\label{recu}
%\noindent
%\vspace{-0.1cm}
%\subsection{Matrix Inversion Lemma for Moving Window}
\label{lem}
\noindent
%Consider the following generalization of the matrix inversion lemma for moving window provided in \cite{lem}.
Consider the following generalization of the matrix inversion lemma for moving window described in \cite[p.~280]{lem}.
\\
\textit{Lemma} \textit{(batch update in moving window)}.~~
Updates of the invertible $n \times n$ matrix $A$, associated with $d$ additions and $r$ removals of the outer products of the vectors $x_j \in \mathbb{R}^n$ and $y_j \in \mathbb{R}^n$ can be expressed in the following equivalent forms:
\begin{align}
   & A = B  + \underbrace{\sum_{j=1}^{d} x_j~x_j^{T}}_{\text{updating}} - \underbrace{\sum_{j=1}^{r}y_j~y_j^{T}}_{\text{downdating}} \label{fform} \\
   & A  = B  + Q~D~Q^{T}   \label{sform}
\end{align}
where $Q = [x_1~x_2~...~x_d~y_1~y_2~\cdots~y_r]$ is the augmented regressor matrix which contains regressor column vectors,
$D = \operatorname{diag}[\underbrace{1,~1~,...~,1}_{\text{d~additions}},\underbrace{-1,-1~,...,~-1}_{\text{r~removals}}]$ is
the signature matrix, which determines the sign structure of the low rank correction, $d + r \ll n$.
\\
Subsequently, the inverse of $A$ can be efficiently obtained through the batch update by applying the Woodbury identity, \cite[p.~350]{bern2} to the compact form (\ref{sform}) as follows:
\begin{align}
   & A^{-1} = B^{-1} - B^{-1}~Q~U^{-1}~Q^{T}~B^{-1}  \label{minv1}
\end{align}
where the update matrix $U = D + Q^{T}~B^{-1}~Q$ is invertible and the inverse $B^{-1}$ is known.
\\ \\
\textit{Corollary.}~~ Rank two updates of the matrix $A$ is a special case of
(\ref{fform}),(\ref{sform}) with $d=r=1$, which can be presented in the form (\ref{sform}) with the following
regressor matrix $Q = [x_1~y_1]$ and the signature matrix $D = \operatorname{diag}[1,-1]$, \cite{sto2023} -  \cite{sto2024}.
% prog lemma1 for different matrcies and sizes
\\ \\
The inverse $A^{-1}$ can also be computed iteratively using the Sherman–Morrison formula, applying it once per rank-one update and using each intermediate inverse in the subsequent step, \cite{lem}, \cite{zhang}. However, this approach comes with a number of practical and theoretical limitations, which include:

\begin{itemize}
    \item  Round-off Error Accumulation: Each application of the Sherman–Morrison formula introduces numerical errors due to finite-precision arithmetic. When applied repeatedly, these errors can accumulate significantly, resulting in a poor approximation of the true inverse.

    \item Loss of Symmetry and Positive Definiteness: Even if the original matrix $A$ is symmetric positive definite (SPD), the intermediate matrices generated during the iterative process may not preserve these properties. This can lead to unstable or invalid intermediate results.

    \item Intermediate Singularity: The Sherman–Morrison formula requires that the intermediate matrices remain invertible at each step. This assumption is not always valid, particularly in the presence of near-linear dependencies in the rank-one updates. Additional assumptions or regularization techniques may be needed to guarantee invertibility.

    \item Poor Cache Performance: The iterative nature of the algorithm leads to frequent and irregular memory accesses, especially when updating and reusing matrix elements stored in memory. This behavior results in poor utilization of CPU cache lines, causing frequent cache misses and slowing down execution.

    \item Repeated Memory Access: Each step of the update involves accessing and modifying elements of the current inverse matrix. Without efficient memory locality, this leads to high memory bandwidth consumption and increased latency due to repeated reads and writes to main memory.

    \item Limited Parallelism: The update steps must be performed sequentially, as each depends on the result of the previous one. This inherent data dependency severely limits opportunities for parallelization, making the method inefficient on modern multicore architectures.

    \item Computational Inefficiency: Due to all the above factors, especially poor cache behavior and lack of parallelism, the iterative Sherman–Morrison approach is often less efficient in practice than direct batch inversion techniques, particularly for large-scale matrices.
\end{itemize}
The inversion method presented in the lemma applies the correction in a single pass, thereby reducing the number of matrix multiplications and mitigating the propagation of numerical errors. It preserves matrix symmetry, improves the conditioning of the resulting inverse, and is inherently more parallelizable. In addition, the method exhibits enhanced memory efficiency, making it particularly effective for deployment in high performance computing architectures.

\section{Recursive Estimation Framework for Segmented Forgetting Profile}
\label{reclr}
\noindent
Assume that the measured oscillating signal  ${y}_k$  and its corresponding model $\hat{y}_k$  are represented as follows:
\begin{align}
 y_k &= \fai_k^{T} \theta_*  + \xi_k  \label{yk} \\
 \hat{y}_k &= \fai_k^{T} \theta_k   \label{yhat}   \\
 \fai_k^{T} &=  [ 1 ~ cos( q_0 k ) ~ sin(q_0 k ) ~ \cdots ~ cos( q_h k ) ~ sin( q_h k ) ] \nonumber
\end{align}
where $\theta_* $ is the vector of unknown parameters, $\theta_k$  is the vector of adjustable parameters, $\fai_k$  is the harmonic regressor, $q_0,...q_h$ are the frequencies, and $\xi_k$ is white Gaussian zero mean noise uncorrelated
with $\fai_k$,  $k = 1,2,...$.
\\ Recursive estimation framework for segmentation of the forgetting profile can be presented as follows:
\begin{align}
& A_k  =  \lambda~A_{k-1} + Q_k~D~Q^{T}_k, ~~ k \ge w + 1 \label{imupdate} \\
& Q_k  = [\fai_{k} ~ \sqrt{|\beta-\lambda|} \fai_{k-1} \cdots
 \sqrt{|\lambda^m - \beta^p| \lambda}  \fai_{k-p-1} ~  \sqrt{\lambda^{m+w-p}} \fai_{k-w}]  \nonumber \\
& D = \operatorname{diag}[ \underbrace{1, ~\operatorname{sign}(\beta-\lambda),~\cdots,}_{p+1}~\operatorname{sign}(\lambda^m - \beta^p),~-1] \nonumber \\
& \Gamma_{k} = \frac{1}{\lambda} ~ [~ \Gamma_{k-1} - \Gamma_{k-1} ~ Q_k ~ S^{-1} ~ Q^T_k ~ \Gamma_{k-1}~], ~~\Gamma_{w} = A^{-1}_w   \label{gk} \\
&  \theta_k  = \theta_{k-1} - \Gamma_{k-1}  ~ Q_k ~ S^{-1} ~ [ Q^T_k~ \theta_{k-1} - \tilde{y}_k ]   \label{tk} \\
&  S = \lambda~D + Q^T_k ~ \Gamma_{k-1} ~ Q_k          \nonumber   \\
&  \tilde{y}^T_k = [y_k ~ \sqrt{|\beta-\lambda|}y_{k-1}~\cdots~ \sqrt{|\lambda^m - \beta^p| \lambda} y_{k-p-1}  ~ \sqrt{\lambda^{m+w-p}} y_{k-w} ] \nonumber
\end{align}
where equation (\ref{imupdate}) represents rank $p+3$ ~ ($p=1,2, \cdots \ll  w$) updates of the information matrix  $A_k$  in the moving window of the size $w$, and the augmented regressor matrix $Q_k$ contains scaled column regressor vectors.
The recursive framework (\ref{imupdate})-(\ref{tk}) supports various types of segmentation  parameterized by  $\beta$, $\lambda$, $m$ and $p$, see Figure~\ref{figseg1}. Each profile yields a distinct sign structure of the signature matrix $D$, which governs
the addition and removal of updates.
\\ The RLS algorithms with low rank updates, which are derived by application of the matrix inversion lemma described in Section~\ref{lem}, can be  written in the form (\ref{gk}), (\ref{tk}), where $\tilde{y}_k$ is the augmented output,  see also \cite{sto2023}, \cite{sto2024}.
\\ The algorithm downweights rapidly recent $p+1$ measured points with the relatively small forgetting factor $ 0 < \beta < 1$, see Figure~\ref{figseg1}. Integer $m > 0 $ determines the magnitude of the drop in the profile so that $\lambda^{m+1} < \beta^p $ and the tail in the segmented profile downweights  slowly older data with larger forgetting factor, $ 0 < \lambda < 1$ for reduction of the condition number of the information matrix.
\\ \\ The details associated with derivation of the recursive estimation framework can be found in the following paper  \cite{lowrank} :
\\ \\
Stotsky A. (2025). Accelerating with low rank updates: RLS estimation with segmentation of the forgetting profile. International Journal of Computer Applications (0975 - 8887), 187(54),1–5.
\\
\url{https://www.ijcaonline.org/archives/volume187/number54/stotsky-2025-ijca-925940.pdf}
\\  \\
The unbiasedness of the algorithm (\ref{gk}), (\ref{tk}),~$\mathbb{E}[\theta_k] = \theta_*$ can be formally derived using the framework presented in \cite{sto2024}, assuming a full-rank information matrix.
The variance of the parameter estimation error is inversely proportional to the eigenvalues of the information matrix. When the information matrix is ill-conditioned, some directions in parameter space which correspond to small eigenvalues lead to large estimation variances along these directions. This typically happens when the regressor signals are not sufficiently excited over the window. In contrast, well-conditioned information matrices result in more uniform and lower variance across all parameters.
\\ Notice that applying an infinite window over a large number of iterations can introduce bias, increase the variance of parameter mismatch, and cause imbalance. These effects often arise due to ill-conditioning, numerical error accumulation, sensitivity to measurement noise, outliers, and unmodeled disturbances.
\\
In contrast, segmented forgetting in the finite window alleviates ill-conditioning and reduces bias, variance, and imbalance in parameter estimation, thereby improving the robustness and performance of the algorithms, see Section~\ref{win}.

\begin{figure*}[!ht]
%\begin{center}
\hspace{-2cm} % or -2cm to move left
  \includegraphics[height=6cm,width=20cm]{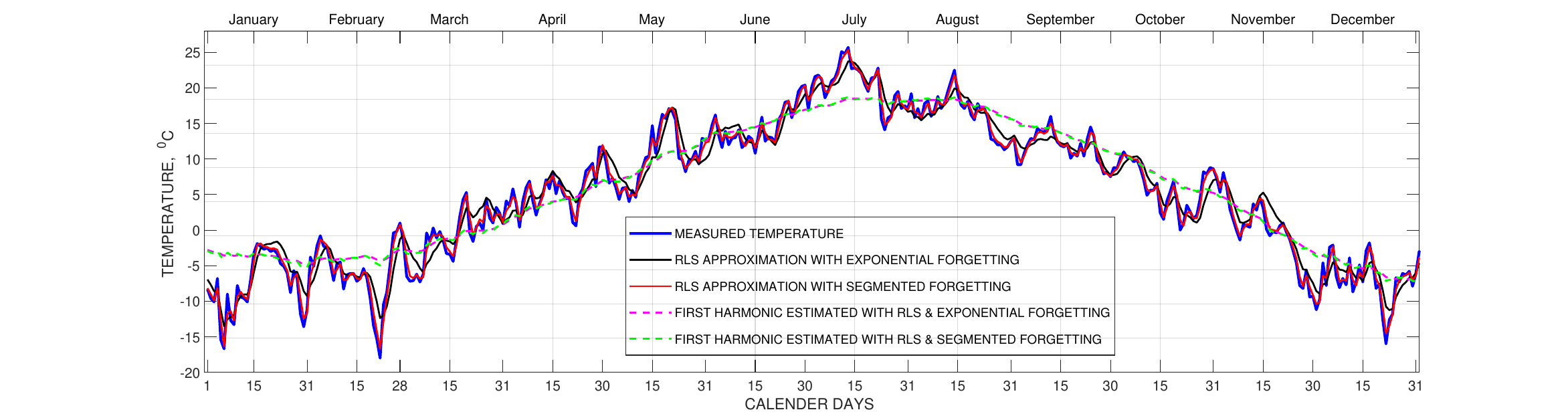}
%  \includegraphics{sto1.pdf}
%\end{center}
\caption{\small{Comparison of  the approximation performance of segmented and exponential forgetting profiles in moving window of the size $w=400$ is shown in this Figure. Daily temperature measurements are plotted with the blue line. The output of the RLS algorithm with rank two updates and $\lambda = 0.99$ is plotted with the black line. The output of RLS algorithm with segmented profile and rank four updates, designed for $p=1$,~$\beta = 0.89$,~$\lambda = 0.99$, $m=250$ (see Figure~\ref{figseg1}) is plotted with the red line. Histograms of approximation errors are presented in Figure~\ref{figm1}.
Estimates of the first harmonics are plotted with magenta and green lines for exponential and segmented profiles respectively.
 }}
\label{figcimer1}
\end{figure*}
\begin{figure*}[!ht]
\begin{center}
%\hspace{-2cm} % or -2cm to move left
  \includegraphics[height=8cm]{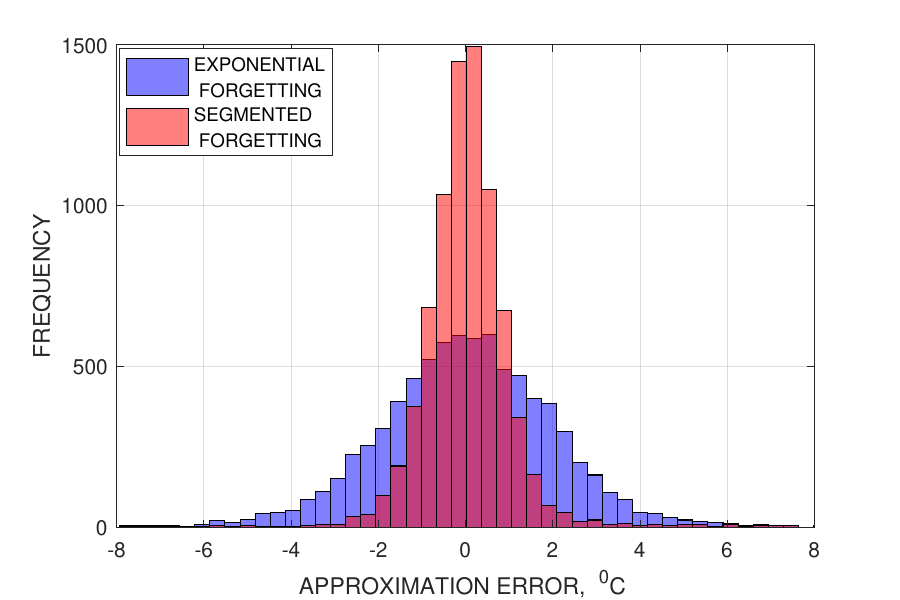}
\end{center}
\caption{\small{Histogram of the approximation error for the exponential forgetting profile is plotted with the blue color, and
the histogram of the approximation error for the segmented forgetting profile is plotted with the red color. Approximation performance is significantly improved via segmentation of the profile.    }}
\label{figm1}
\end{figure*}

\section{Temperature Forecasts}
\label{int}
\begin{figure*}[!ht]
%\begin{center}
\hspace{-2cm} % or -2cm to move left
  \includegraphics[height=6cm,width=20cm]{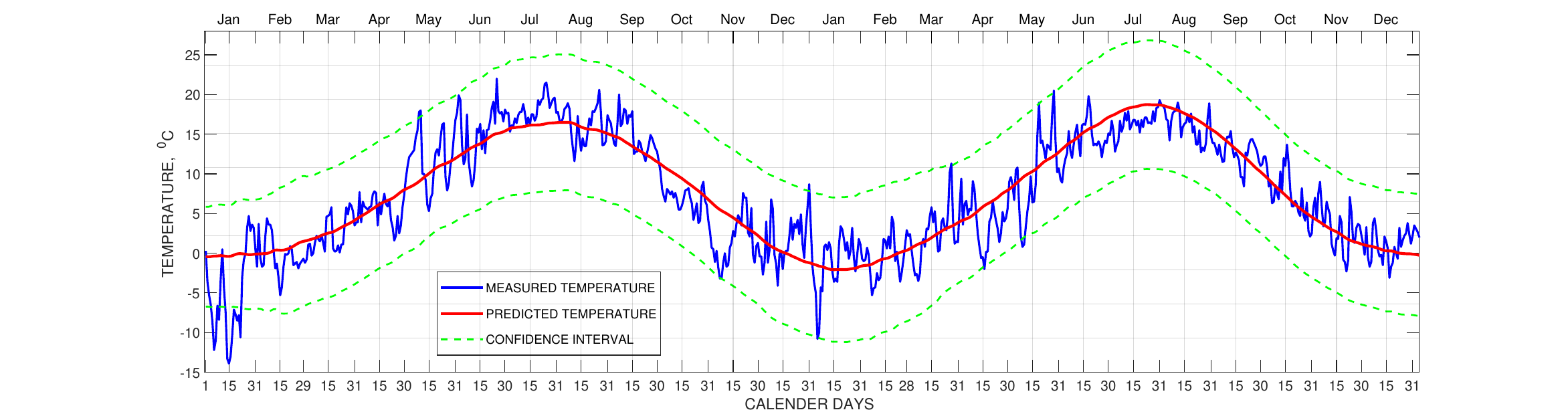}
%  \includegraphics{sto1.pdf}
%\end{center}
\caption{\small{The Figure shows the $30$ days temperature forecast based on prediction of the first harmonic and three sigma confidence interval with estimation of the variance in moving window. It is shown that the seasonal trend can be very well predicted using estimate of the first harmonic of low resolution daily temperature measurements.
The prediction estimates the mean, maximum, and minimum temperature values projected 30 days ahead.
The accuracy of the prediction is assessed by checking if the actual temperature measurements lie within the confidence intervals established around the predicted first harmonic values. Most observed temperature measurements fall within the confidence intervals of the first harmonic predictions, confirming their reliability.
 }}
\label{figcimer2}
\end{figure*}

\noindent
Temperature forecasts covering a period longer than one week are essential for providing early insights into weather trends that facilitate effective planning and preparedness. Such forecasts contribute to improved decision-making across multiple sectors, including agriculture, energy, transportation, and public health, by offering a comprehensive outlook on forthcoming climate conditions. For instance, farmers rely on extended temperature forecasts to schedule planting and harvesting, energy providers use them to anticipate demand fluctuations, and public health agencies use them to prepare for heatwaves or cold snaps.
\\
One of the key challenges in producing accurate long-range temperature forecasts lies in capturing the inherent variability of the atmosphere while also recognizing the underlying patterns that influence it. This is where the concept of periodicity becomes critically important. Temperature naturally follows repeating patterns, such as daily, weekly, and seasonal cycles driven by predictable astronomical and environmental factors. Daily cycles are shaped by the Earth's rotation, while seasonal cycles arise from its orbit around the sun and the resulting variation in solar radiation. Weekly patterns, though less rooted in natural law, can reflect consistent human activities that affect urban heat generation and energy consumption.
\\
Recognizing and modeling these periodic components allows for more robust predictions. System identification techniques, mathematical tools used to model dynamic systems based on input and output data are particularly effective in detecting these recurring patterns within temperature time series data. By applying these techniques, forecasters can decompose temperature signals into their periodic and aperiodic components, allowing for clearer interpretation and projection of future behavior.
\\
Leveraging periodicity in this manner significantly enhances forecasting accuracy by filtering out noise and isolating consistent trends. As a result, long-term temperature predictions become not only more reliable but also more actionable, enabling stakeholders to make informed decisions well in advance. Ultimately, integrating periodic analysis through system identification into forecasting frameworks represents a vital step forward in the evolution of climate prediction capabilities.
\\
Assume that the temperature time series is represented by equation (\ref{yk}), incorporating the fundamental frequency associated with the annual cycle, as well as sixteen higher order harmonics that capture additional periodic components of shorter durations.
The model of the signal is presented in the form (\ref{yhat}) with adjustable parameters  (\ref{gk}),(\ref{tk}).
To validate the model, low resolution daily mean temperature measurements from the Stockholm Old Astronomical Observatory, \cite{data1} are used.
\\ Comparison of the approximation performance of segmented and exponential forgetting profiles is shown in Figure~\ref{figcimer1} and  Figure~\ref{figm1}.
%Measured temperature, which is plotted with the blue line is approximated by RLS algorithm with rank two updates, whose output is plotted with the black line and by RLS algorithm with segmented forgetting profile and rank four updates, whose output is plotted with the red line.
The measured temperature, shown by the blue line in Figure~\ref{figcimer1}, is approximated using two RLS algorithms: one with rank two updates (black line), and another with segmented forgetting and rank four updates (red line).
The histograms in Figure~\ref{figm1} illustrate that the approximation performance is significantly improved by segmenting the forgetting profile.
%Moreover, the first harmonic which corresponds to the annual periodicity is estimated more accurately, which has a direct impact of the long term temperature forecast, presented in Figure~\ref{figcimer2}.
Moreover, the first harmonic, which corresponds to the annual periodicity, is estimated more accurately, directly impacting the long term temperature forecast presented in Figure~\ref{figcimer2}.
\\ The Figure~\ref{figcimer2}  presents the 30-day-ahead temperature forecast based on the first harmonic component, accompanied by a three sigma confidence interval, \cite{card}.  The variance is also estimated within the moving window.
It has been demonstrated that the seasonal trend can be accurately predicted using estimates of the first harmonic obtained from low resolution temperature measurements.
%Almost all actual temperature measurements fall within the established confidence intervals around predicted values of the first harmonic that confirms reliability of the predictions.
Nearly all observed temperature measurements fall within the established confidence intervals around the predicted values of the first harmonic, confirming the reliability of the predictions.
%The accuracy of the prediction is evaluated by verifying whether the actual temperature measurements fall within the established confidence intervals around predicted values of the first harmonic.

\section{Conclusion \& Outlook}
\label{con}
\noindent
The regularization method, aligned with effective segmentation of the forgetting profile, showed strong potential in optimizing the balance between rapidity, desired condition number of the information matrix, accuracy, and numerical stability.
The main contribution of this work is the integration of the segmented profile into a low rank update recursive least squares framework. The development utilizes the matrix inversion lemma tailored for moving window computations.
\\
Design flexibility in the forgetting profile enabled more accurate approximation of the frequency content in low resolution temperature measurements, thereby improving the reliability of temperature predictions.
\\
Further investigation of the properties of RLS algorithms with low rank updates is warranted to achieve additional performance improvements, in accordance with the concept of the profile segmentation introduced in this report.
\\
Finally, the recursive framework with low rank updates
facilitates the development and assessment of new, more advanced segmentation strategies which are anticipated to further improve estimation performance. Moreover, performance can be further enhanced through the use of Newton-Schulz and Richardson corrections,\cite{sto2023}, \cite{ifac2022}.

\section {Disclosure Statement}
This research was not supported by any organization.

\end{document}